\documentclass[12pt,draft]{amsart}
\usepackage{amsfonts,amsmath,amsthm,latexsym,amscd,amsbsy,amssymb,fleqn,leqno,pb-diagram,lamsarrow}

\chardef\bslash=`\\ 

\makeatletter
\def\verbatim{\interlinepenalty\@M \@verbatim
  \leftskip\@totalleftmargin\advance\leftskip2pc
  \frenchspacing\@vobeyspaces \@xverbatim}
\makeatother
\makeatletter
  \def\dgt@k{\dg@DX=-3 \dg@DY=2 \dg@SIZE=3}
\makeatother

\makeatletter
  \def\dgt@kk{\dg@DX=3 \dg@DY=-1 \dg@SIZE=3}%
\makeatother

\theoremstyle{plain}
\newtheorem{thm}{Theorem}[section]
\newtheorem{cor}[thm]{Corollary}
\newtheorem{lem}[thm]{Lemma}
\newtheorem{pro}[thm]{Proposition}

\theoremstyle{definition}

\numberwithin{equation}{section}

\newcounter{rmnum}

\def\symbolnote#1#2{\let\thefootn=\thefootnote%
\renewcommand{\thefootnote}{\fnsymbol{footnote}}%
\footnotemark[#1]%
\footnotetext[#1]{#2}%
\let\thefootnote=\thefootn
}

\newfont{\bbb}{msbm10 scaled \magstep1}
\newfont{\bbc}{msbm8 scaled \magstep0}

\newcommand{\uin}{\mbox{\bbb I}}

\newcommand{\ds}{\displaystyle}

\begin{document}


\title{On commutative and non-commutative $C^{\ast}$-algebras with the approximate $n$-th root property}

\author{A. Chigogidze}
\address{Department of Mathematical Sciences,
University of North Carolina at Greensboro,
P.O. Box 26170, Greensboro, NC 27402-6170, U.S.A.}
\email{chigogidze@uncg.edu}

\author{A. Karasev}
\address{Department of Computer Science and Mathematics,
Nipissing University,
100 College Drive, P.O. Box 5002, North Bay, ON, P1B 8L7, Canada}
\email{alexandk@nipissingu.ca}
\thanks{The second author was partially supported by his NSERC Grant 257231-04.}

\author{K. Kawamura}
\address{Institute of Mathematics,
University of Tsukuba, Tsukuba, Ibaraki 305-8071, Japan}
\email{kawamura@math.tsukuba.as.jp}
\thanks{The paper was started during the third author's visit to Nipissing University in July 2004}

\author{V.  Valov}
\address{Department of Computer Science and Mathematics, Nipissing University,
100 College Drive, P.O. Box 5002, North Bay, ON, P1B 8L7, Canada}
\email{veskov@nipissingu.ca}
\thanks{The last author was partially supported by his NSERC Grant 261914-03}

\keywords{$C^*$-algebras with the approximate $n$-th root property, rank of $C^*$-algebras, covering dimension, invertible maps}
\subjclass{Primary: 46L85; Secondary: 54C40.}


\begin{abstract}
We say that a $C^*$-algebra $X$ has the approximate $n$-th root property 
($n\geq 2$) if for every $a\in X$ with $\|a\|\leq 1$ and every $\varepsilon>0$ 
there exits $b\in X$ such that
$\|b\|\leq 1$ and $\|a-b^n\|<\varepsilon$. 
Some properties of commutative and non-commutative $C^*$-algebras having the 
approximate $n$-th root property are investigated. In particular, it is 
shown that there exists a non-commutative (resp., commutative) separable 
unital $C^*$-algebra $X$ such that any other (commutative) separable unital 
$C^*$-algebra is a quotient of $X$.  Also we illustrate a commutative 
$C^*$-algebra, each 
element of which has a square root such that its maximal ideal space has 
infinitely generated first {\v C}ech cohomology.   
\end{abstract}

\maketitle
\markboth{A. Chigogidze, A. Karasev, K.~Kawamura and V.~Valov}{approximate $n$-th root property}


\section{Introduction}
All topological spaces in this paper are assumed to be (at least) completely 
regular. A compact Hausdorff space is called a {\it compatum} for simplicity.  
By  $C^{\ast}$-algebra and homomorphisms between $C^{\ast}$-algebras, we mean 
unital $C^{\ast}$-algebras and unital $*$-homomorphisms. 
For a space $X$ and an integer $n\geq 2$, we consider the following conditions
($\|\cdot\|$ denotes the supremum norm):

\smallskip\noindent
\begin{itemize}
\item [$\displaystyle (*)_n$] For each bounded continuous function 
$f\colon X\to\mathbb C$ and each $\varepsilon>0$, there exists a continuous 
function $g\colon X\to\mathbb C$ such that $\|f-g^n\|<\varepsilon$.

\item[$\displaystyle (**)_n$] For each bounded continuous function 
$f\colon X\to\mathbb C$ and each $\varepsilon>0$, there exist bounded 
continuous functions $g_1,..,g_n\colon X\to\mathbb C$ such that
$f=\prod_{i=1}^{i=n}g_i$ and $\displaystyle \|g_i-g_j\|<\varepsilon$ for 
each $i,j$.
\end{itemize}

We say that the space $C^{\ast}(X)$ of all bounded complex-valued functions on $X$ has the approximate $n$-th root property if $X$ satisfies condition $(*)_n$. 
The results in this paper were inspired by the following theorem established by K. Kawamura and T. Miura \cite{km}:

\begin{thm}
Let $X$ be a compactum with $\dim X\leq 1$ and $n$ a positive integer. Then the following conditions are equivalent:

\begin{itemize}
\item[(1)] $C(X)$ has the approximate $n$-th root property. 
\item[(2)] $X$ satisfies condition $(**)_n$. 
\item[(3)] the first
\v{C}ech cohomology $\displaystyle\check H^1(X;\mathbb Z)$ is
$n$-divisible, that is, each element of
$\displaystyle\check{H}^1(X;\mathbb Z)$ is divided by $n$.
\end{itemize}
\end{thm}

Let $\mathcal{A}(n)$  denote the class of all completely regular spaces 
satisfying condition $(*)_n$ and $\mathcal{A}_1(n)$ is the subclass of $\mathcal{A}(n)$ consisting of spaces $X$ with $\dim X\leq 1$.  Let also 
$\mathcal{H}(n)$  denote the class of all {\it compacta} $X$ with
$\displaystyle\check H^1(X;\mathbb Z)$ being $n$-divisible.  

In Section 2 we investigate some properties of the classes $\mathcal{A}(n)$, $\mathcal{A}_1(n)$ and $\mathcal{H}(n)$. In particular, the following theorem is established:  

\begin{thm}
Let $n$ be a positive integer and let $\mathcal{K}$ denote one of the classes 
$\mathcal{A}(n)$, $\mathcal{A}_1(n)$ or $\mathcal{H}(n)$. Then, for every 
cardinal $\tau\geq\omega$, there exists a compactum $X_{\tau} \in\mathcal{K}$ 
of weight $\leq\tau$ and a $\mathcal{K}$-invertible map 
$f_{\mathcal{K}}\colon X_{\tau} \to\uin^{\tau}$.
\end{thm}

Here, a map $h\colon X\to Y$ is said to be {\it invertible} for the class 
$\mathcal{K}$ (or simply, $\mathcal{K}$-{\it invertible}) if for every map 
$g\colon Z\to Y$ with $Z\in\mathcal{K}$ there exists a map $\overline{g}\colon Z\to X$ such that $g=h\circ\overline{g}$. 

Theorem 1.2 implies the next corollary.

\begin{cor}
Let $n$ be a positive integer and let $\mathcal{K}$ be one of the classes 
$\mathcal{A}(n)$, $\mathcal{A}_1(n)$ or $\mathcal{H}(n)$. Then, for every $\tau\geq\omega$, there exists a compactum $X\in\mathcal{K}$ of weight $\tau$ which contains every space from $\mathcal{K}$ of weight $\leq\tau$.
\end{cor}

It is easily seen that the modification of condition $(*)_n$, obtained by 
requiring both $f$ and $g$ to be of norm $\leq 1$, is equivalent to $(*)_n$. 
This 
observation leads us to consider the following classes of general 
(non-commutative) $C^{\ast}$-algebras.
We say that a $C^{\ast}$-algebra $X$ satisfies {\it the approximation} 
$n$-{\it th root property} if for every $a\in X$ with $\|a\|\leq 1$ and every 
$\varepsilon>0$ there exists $b\in X$ such that
$\|b\|\leq 1$ and $\|a-b^n\|<\varepsilon$. The class of all 
$C^{\ast}$-algebras with the approximate $n$-th root property is denoted by 
$\mathcal{AP}(n)$. Let $\mathcal{AP}_1(n)$ be the subclass of
$\mathcal{AP}(n)$ consisting of $C^{\ast}$-algebras of bounded rank $\leq 1$ 
(recall that  bounded rank of $C^*$-algebras is a non-commutative analogue of 
the covering dimension $\dim$, see \cite{cv}). We also consider the class 
$\mathcal{HP}(n)$ of $C^{\ast}$-algebras $X$ with the following property: for 
every {\it invertible} element $a\in X$ with $\|a\|\leq 1$ and every 
$\varepsilon>0$ there exists 
$b\in X$ such that $\|b\|\leq 1$ and $\|a-b^n\|<\varepsilon$.

\medskip\noindent
In the sequel, $\mathcal{AP}(n)_s$ denotes the class of all separable 
$C^*$-algebras from $\mathcal{AP}(n)$. 
The notations $\mathcal{AP}_1(n)_s$ and $\mathcal{HP}(n)_s$ have the same 
meaning.

\medskip\noindent
Recall now the concept of $\mathcal{\Re}$-invertibility introduced
in \cite{ch3}, where $\mathcal{\Re}$ is a given class of
$C^{\ast}$-algebras. A homomorphism $p\colon X\to Y$ is said to be 
$\mathcal{\Re}$-{\it invertible} if, for any homomorphism
$g\colon X\to Z$ with $Z\in\mathcal{\Re}$, there exists a
homomorphism $\overline{g}\colon Y\to Z$ such that
$g=\overline{g}\circ p$. We also introduce the notion of a {\it universal} 
$C^*$-algebra for a given class $\mathcal{\Re}$ as a $C^*$-algebra 
$Y\in\mathcal{\Re}$ such that any other $C^*$-algebra from $\mathcal{\Re}$ 
is a quotient of $Y$. 

Section 3 is devoted to the classes $\mathcal{AP}(n)$, $\mathcal{AP}_1(n)$ and 
$\mathcal{HP}(n)$. The results of this section can be considered as 
non-commutative counterparts of the results from Section 2. For example, 
Theorem 1.4 below is a non-commutative version of Theorem 1.2.

\begin{thm} Let $n$ be a positive integer and let $\mathcal K$ be one of the 
classes $\mathcal{AP}(n)$, $\mathcal{AP}_1(n)$ and $\mathcal{HP}(n)$. Then
there exists a ${\mathcal K}$-invertible
unital $\ast$-homomorphism
$p \colon C^{\ast}\left( {\mathbb F}_{\infty}\right)
\to Z_{\mathcal K}$ of $C^{\ast}(F_{\infty})$ to a separable unital $C^{\ast}$-
algebra $Z_{\mathcal K} \in {\mathcal K}$, where $C^{\ast}({\mathbb F}_{\infty})$
is the group 
$C^{\ast}$-algebra of 
the free group on countable number of generators.
\end{thm}

It is well-known that every separable
$C^{\ast}$-algebra is a surjective image of $C^{\ast}({\mathbb F}_{\infty})$. 
Therefore, if  $\mathcal{\Re}$ is a class of separable $C^*$-algebras and 
$p\colon C^{\ast}(F_{\infty})\to Y_{\mathcal{\Re}}$ is a
$\Re$-invertible homomorphism with $Y_{\mathcal{\Re}}\in\mathcal{\Re}$, then
$Y_{\mathcal{\Re}}$ is universal for the class $\mathcal{\Re}$. Hence, 
Theorem 1.4 implies 
that each of the classes $\displaystyle\mathcal{AP}(n)_s$, 
$\mathcal{AP}_1(n)_s$ and $\mathcal{HP}(n)_s$ has a universal element.

Let us note that there exists a non-commutative separable $C^*$-algebra which 
belongs to any one of the classes $\displaystyle\mathcal{AP}(n)$, 
$\mathcal{AP}_1(n)$ and $\mathcal{HP}(n)$. Indeed, let $X=M(m)$ be the algebra 
of all $m\times m$ complex matrixes, where $m\geq 2$ is a fixed integer.   By 
\cite{bp}, the bounded rank of $X$ is 0. Moreover,  using the canonical Jordan 
form representation, one can show that if $A\in X$ and $n\geq 2$, then
$A$ can be approximated by a matrix $B\in X$ with  $C^n=B$ for some $C\in X$. 
Hence,  
$X$ is a common element of $\displaystyle\mathcal{AP}(n)$, $\mathcal{AP}_1(n)$ 
and $\mathcal{HP}(n)$. This implies that the universal elements of 
$\displaystyle\mathcal{AP}(n)_s$, $\mathcal{AP}_1(n)_s$ and 
$\mathcal{HP}(n)_s$ are also non-commutative.     

Section 4 deals with {\it square root closed compacta}, compacta 
$X$ such that, for every $f\in C(X)$, there is $g\in C(X)$ with $f=g^2$. 
It is known that if $X$ is a first-countable connected compactum, then $X$ is 
square-root closed if and only if $X$ is locally connected, $\dim X\leq 1$ and 
$\displaystyle\check{H}^1(X;\mathbb Z)$ is trivial, see \cite{cm}, 
\cite{hm}, \cite{km} and \cite{mn}. A topological characterization of 
general square root closed compacta has not been known. 
Here we show that a square root closed compactum $X$ with $\dim X\leq 2$, 
constructed based on the idea of B. Cole (\cite{stout}, Chap.3, section 19) 
and M.I. Karahanjan \cite{kar} has infinitely generated first 
{\v C}ech cohomology 
$\displaystyle\check{H}^1(X;\mathbb Z)$. This space $X$ is the limit space of 
an inverse system $\big(X_{\alpha}, \pi^{\beta}_{\alpha}:\alpha<\omega_1\big)$ 
starting with the unit disk in the plane and such that each map 
$\pi^{\beta}_{\alpha}\colon X_{\beta}\to X_{\alpha}$ is invertible with 
respect to the class of square root closed compacta. A similar 
construction yields a one-dimensional such compactum.  This illustrates that 
the topological characterization of (not necessarily first countable) square 
root closed compacta would be rather different 
than the one for first-countable compacta mentioned above.  Also, the 
invertibility $\pi^{\beta}_{\alpha}\colon X_{\beta}\to X_{\alpha}$  
allows us to obtain a universal element for the class of square root closed 
compacta with arbitrarily fixed weight. 

\section{Some properties  of the classes $\mathcal{A}(n)$, $\mathcal{A}_1(n)$ and $\mathcal{H}(n)$}

\begin{lem} For a compactum $X$, the following conditions are equivalent:

\begin{itemize}
\item[(1)] For any $f\colon X\to S^1$ and any $\varepsilon>0$ there exists $g\colon X\to S^1$ such that $\|f-g^n\|<\varepsilon$.
\item[(2)] $\displaystyle\check{H}^1(X;\mathbb Z)$ is $n$-divisible.
\end{itemize}
\end{lem}

\begin{proof}
When $\varepsilon=0$ in (1), this equivalence was established by 
Kawamura-Miura in \cite[Lemma 3.1]{km}. Their arguments remain also valid in 
the present situation because any two sufficiently close functions from $X$ 
into $S^1$ are homotopic.
\end{proof}

\begin{lem}
Let $X$ be the limit space of an inverse system 
$\{X_{\alpha},p^{\beta}_{\alpha}: \alpha, \beta\in A\}$ of compacta. 
Then, for every $f\in C(X)$ and every $\varepsilon>0$, there exists 
$\alpha\in A$ and $g\in C(X_\alpha)$ such that $g\circ p_\alpha$ is 
$\varepsilon$-close to $f$, where $p_\alpha\colon X\to X_\alpha$ is the 
$\alpha$-th limit projection. Moreover, $g\in C(X,S^1)$ provided
$f\in C(X,S^1)$.
\end{lem}

\begin{proof}
We take a finite cover $\omega$ of $f(X)$ consisting of open and convex 
subsets of $\mathbb C$ each of diameter $<\varepsilon$. Since $X$ is compact, 
we can find $\alpha$ and an open cover $\gamma=\{U_j:j=1,..,m\}$ of $X_\alpha$ 
such that $p_{\alpha}^{-1}(\gamma)$ is a star-refinement of the cover 
$f^{-1}(\omega)$. Without loss of generality, we can assume that each $U_j$ 
is functionally open in $X_\alpha$, i.e., $U_j=h_j^{-1}((0,1])$ for some 
function $h_j\colon X_{\alpha}\to [0,1]$. For any $j$ we fix a point 
$x_j\in p_\alpha^{-1}(U_j)$ and the required function 
$g\colon X_\alpha\to\mathbb C$ is defined by 
$\displaystyle g(y)=\sum_{j=1}^{j=m}h_j(y)f(x_j)$. When $f\in C(X,S^1)$ and  
$\varepsilon$ is sufficiently small, 
$g(X_\alpha)\subset\mathbb C\backslash\{0\}$  and, by considering the 
composition of $g$ and the usual retraction 
$r\colon\mathbb C\backslash\{0\}\to S^1$, we can assume $g\in C(X_\alpha,S^1)$.
\end{proof}

\begin{cor}
Let $\mathcal K$ be one of the classes $\mathcal{A}(n)$, $\mathcal{A}_1(n)$ 
and $\mathcal{H}(n)$. If
$X$ is the limit space of an inverse system 
$\{X_{\alpha},p^{\beta}_{\alpha}: \alpha, \beta\in A\}$ of compacta with each 
$X_\alpha\in\mathcal K$, then $X\in\mathcal K$.
\end{cor}

\begin{proof}
This is a direct application of Lemma 2.2 for the class $\mathcal{A}(n)$. 
Since the limit space of any inverse system of at most one dimensional 
compacta is of dimension $\leq 1$, the validity of our corollary for 
$\mathcal{A}(n)$ yields its validity for $\mathcal{A}_1(n)$. Finally, 
Lemma 2.1 and Lemma 2.2 settle the proof for the class  $\mathcal{H}(n)$.
\end{proof}

We say that a class of spaces $\mathcal{K}$ is {\it factorizable} if, for 
every map $f\colon X\to Y$ of a {\it compactum} $X \in \mathcal{K}$, there 
exists a compactum $Z\in\mathcal{K}$ of weight $w(Z)\leq w(Y)$ and maps 
$\pi\colon X\to Z$ and
$p\colon Z\to Y$ such that $f=p\circ\pi$.

\begin{pro}
Any one of the classes $\mathcal{A}(n)$, $\mathcal{A}_1(n)$ and 
$\mathcal{H}(n)$ is factorizable.
\end{pro}

\begin{proof} We consider first the class $\mathcal{H}(n)$.  Fix a map
$\displaystyle f\colon X\to Y$ of a compactum $X\in\mathcal{H}(n)$ and 
assume $w(Y)\leq\tau$. Obviously, we can assume $X$ is of weight $w(X)>\tau$ 
and $Y$ is compact. By induction, we construct sequences of compacta $X_k$, 
dense subsets $\displaystyle M_k\subset C(X_k,S^1)$ of cardinality $\leq\tau$ 
and maps $\displaystyle\pi_k\colon X\to X_k$,
$\displaystyle p^{k+1}_k\colon X_{k+1}\to X_k$, $k\geq 0$, satisfying the 
following conditions:

\begin{itemize}
\item[(0)] $X_0=Y$, $\pi_0=f$, 
\item[(1)] $p^{k+1}_k\circ\pi_{k+1}=\pi_k$, $w(X_k)\leq\tau$ and $M_{k}$ 
separates points of $X_k$ ($k\geq 0$);
\end{itemize}

\begin{itemize}
\item[(2)] For every $h\in M_k$ and every $\varepsilon>0$, there exists 
$g\in M_{k+1}$ such that $\|h\circ p^{k+1}_k-g^n\|<\varepsilon$ ($k \geq 0$).
\end{itemize}

\bigskip
The weight of the function space $C(Y,S^1)$ is $\leq\tau$, so $C(Y,S^1)$ 
contains a dense subset $M_0$ of cardinality $\leq\tau$, separating points of 
$Y$.
Suppose the spaces $X_i$, the sets $M_i$ and the maps $\pi_i$, 
$p^i_{i-1}$, $i\leq k$, have been constructed for some $k$.
Since $X\in\mathcal{H}(n)$,
for each $h\in M_k$ and each positive rational number $r\in Q^{+}$, there 
exists $\displaystyle g(h,r)\in C(X,S^1)$ with
$\displaystyle \|h\circ\pi_k-g(h,r)^n\|<r$. Let 
$\pi_{k+1} \colon X \to X_{k} \times (S^{1})^{M_{k} \times Q^{+}} \times 
(S^{1})^{M_{k}}$ 
be the diagonal product of $\pi_k$ and all maps $g(h,r)$ and $h\circ\pi_k$, 
where $h\in M_k$, $r\in Q^{+}$. Let $X_{k+1}=\pi_{k+1}(X)$ and 
$p^{k+1}_k \colon X_{k+1} \to X_{k}$ be the natural projection onto $X_k$.
Since $M_k$ separates points of $X_k$ (condition (1)), $\pi _{k+1}$ is an 
embedding and hence every $g(h,r)$ can 
be represented as $g_{k+1}(h,r)\circ\pi_{k+1}$ with 
$g_{k+1}(h,r)\in C(X_{k+1},S^1)$.
Because $w(X_{k+1})\leq\tau$, $C(X_{k+1},S^1)$ contains a dense subset 
$M_{k+1}$ of cardinality $\leq\tau$ containing all
$g_{k+1}(h,r)$,  $h\in M_k$, $r\in Q^+$ and also separating points of 
$X_{k+1}$. 
Obviously, $X_{k+1}$, $M_{k+1}$ and $\pi_{k+1}$ satisfy conditions 
(1) and (2). Let $Z$ be the limit of the inverse sequence 
$\{X_k,p^{k+1}_k: k=1,2..\}$, $p\colon Z\to Y$ the  first limit projection and 
$\pi\colon X\to Z$ the limit of the maps $\pi_k$. Also let $p_{k} \colon 
Z \to X_{k}$ be the $k$-th limit projection. By Lemma 2.2, for every 
$h\in C(Z,S^1)$ and every $\varepsilon>0$, there exists $m$ and 
$g_m\colon X_m\to S^1$ such that $\|h-g_m\circ p_m\|<\varepsilon/3$.  Now, 
take $h_m\in M_m$ with $\|g_m-h_m\|<\varepsilon/3$. According to our 
construction, $\|h_m\circ p^{m+1}_m-g^n\|<\varepsilon/3$ for some 
$g\in M_{m+1}$. Hence, $\|h-(g\circ p_{m+1})^n\|<\varepsilon$. Finally, 
by Lemma 2.1, we see $Z\in\mathcal{H}(n)$.

The same arguments remain valid when the class $\mathcal{H}(n)$ is replaced by 
$\mathcal{A}(n)$. The only difference is that we have to consider the 
function spaces $C(X_k)$ instead of $C(X_k,S^1)$.  For the class 
$\mathcal{A}_1(n)$ we need the following modifications: all $M_k$, $k\geq 0$, 
are dense subsets of $C(X_k)$ of cardinality $|M_k|\leq\tau$ satisfying 
conditions (1) and (2), where the compactum $X_k$ is of dimension $\leq 1$ for 
each $k \geq 1$. It suffices to demonstrate the construction of $X_1$ and 
$M_1$. Using the above notations, take the diagonal product 
$q_{1} \colon X \to Y \times {\mathbb C}^{M_{0} \times Q^{+}} \times 
{\mathbb C}^{M_0}$
 of $\pi_0=f$ and all maps $g(h,r)$ and $h\circ\pi_0$, where $h\in M_0$ and 
$r\in Q^+$. Let also $Z_1=q_1(X)$ and $q_0\colon Z_1\to Y$ be the natural 
projection. Then, $w(Z_1)\leq\tau$ and, by the Marde\v{s}i\v{c} factorization 
theorem \cite{m}, there exists a compactum $X_1$ of weight $\leq\tau$ and 
$\dim X_{1} \leq 1$, and maps $\pi_1\colon X\to X_1$ and 
$q_2\colon X_1\to Z_1$ with $q_1=q_2\circ\pi_1$. Obviously, every 
$g(h,r)$ can be represented as $g_1(h,r)\circ\pi_1$ with $g_1(h,r)\in C(X_1)$. 
We denote $p^1_0=q_0\circ q_2$ and choose a dense subset $M_1\subset C(X_1)$ 
such that $|M_1|\leq\tau$ and $M_1$ contains every $g_1(h,r)$ with $h\in M_0$ 
and $r\in Q^+$, and separates points of $X_1$. In this way we obtain the 
spaces $X_k$ with $\dim X_k\leq 1$. The last inequalities imply that the limit 
space $Z$ is also of dimension $\leq 1$. Moreover, by Lemma 2.2,
$Z$ satisfies $(*)_n$, so $Z\in\mathcal{A}_1(n)$.
\end{proof}

\begin{cor}
Let $\mathcal K$ be one of the classes $\mathcal{A}(n)$ and 
$\mathcal{A}_1(n)$. Then every space $X\in\mathcal K$ has a compactification 
$Z\in\mathcal K$ with $w(Z)=w(X)$.
\end{cor}

\begin{proof}
Obviously, $X\in\mathcal K$ implies $\beta X\in\mathcal K$. Let $Y$ be an 
arbitrary compactification of $X$ with $w(Y)=w(X)$ and let 
$f\colon\beta X\to Y$ 
the extension of the identity on $X$. Then, by Proposition 2.4, there exists a 
compactum $Z\in\mathcal K$  and maps $g\colon \beta X\to Z$ and 
$h\colon Z\to Y$ with $h\circ g=f$ and $w(Z)=w(X)$. It remains only to observe 
that $Z$ is a compactification of $X$.
\end{proof}

\begin{pro}
Let $\mathcal K$ be one of the classes $\mathcal{A}(n)$, $\mathcal{A}_1(n)$ 
and $\mathcal{H}(n)$. Then every compactum $X\in\mathcal K$ can be represented 
as the limit space of an $\omega$-spectrum 
$\{X_{\alpha},p^{\beta}_{\alpha}: \alpha, \beta\in A\}$ of metrizable compacta 
with each $X_\alpha\in\mathcal K$.
\end{pro}

\begin{proof} Because of similarity of the arguments, we consider only the 
class $\mathcal{A}(n)$. First,
represent $X$ as the limit space of an $\omega$-spectrum
$\{X_{\alpha},p^{\beta}_{\alpha}: \alpha, \beta\in\Lambda\}$ and introduce 
the relation $L$ on $\Lambda ^2$ consisting of all 
$(\alpha,\beta)\in\Lambda^2$ such that $\alpha\leq\beta$ and for each 
$f\in C(X_\alpha)$ and $\varepsilon>0$ there is $g\in C(X_\beta)$ with 
$\|f\circ p^{\beta}_{\alpha}-g^n\|<\varepsilon$. The relation $L$ has the 
following properties:\\ 
(i) for every $\alpha\in\Lambda$ there exists $\beta\in\Lambda$ with 
$(\alpha,\beta)\in L$:\\ 
(ii)if $(\alpha,\beta)\in L$ and $\beta\leq\gamma$, 
then $(\alpha,\gamma)\in L$;\\ 
(iii) if $\{\alpha_k\}$ is a chain in $\Lambda$ with each 
$(\alpha_k,\beta)\in L$, then $(\alpha,\beta)\in L$, where 
$\alpha=\sup\{\alpha_k\}$.\\ 

Indeed, to show (i), we take a countable dense subset 
$M_\alpha\subset C(X_\alpha)$ and, as in Proposition 2.4, for every 
$h\in M_\alpha$ and $r\in Q^+$ choose
$g(h,r)\in C(X)$ with $\|h\circ p_\alpha-g(h,r)^n\|<r$. Notice that, for each 
$f\in C(X)$, there is a $\gamma\in\Lambda$ and $\varphi\in C(X_\gamma)$ such 
that
$f=\varphi\circ p_\gamma$.  Applying this to $g(h,r)$, we can find 
$\beta\in\Lambda$, $\beta > \alpha$, such that for each 
$(h,r) \in M_{\alpha} \times Q^{+}$, we have 
$g(h,r) = g_{\beta}(h,r) \circ p_\beta$, where 
$g_{\beta}(h,r) \in C(X_\beta)$. Then $(\alpha,\beta)\in L$. Property (ii) 
follows directly and (iii) follows from Lemma 2.2 and the fact that
$X_\alpha$ is the limit space of the inverse sequence generated by 
$\displaystyle X_{\alpha _k}$ and the projections
$p_{\alpha_{k}}^{\alpha_{k+1}} \colon X_{\alpha _{k+1}} \to X_{\alpha _{k}}$, 
$k=1,..$, because  $\alpha$ is supremum of the chain $\{\alpha_k\}$.\\

By \cite[Proposition 1.1.29]{ch1}, the set 
$A=\{\alpha\in\Lambda:(\alpha,\alpha)\in L\}$ is cofinal and 
$\omega$-closed in $\Lambda$. Obviously, $X_\alpha\in\mathcal{A}(n)$ for each 
$\alpha\in A$ and $X$ is the limit of the inverse system 
$\{X_{\alpha},p^{\beta}_{\alpha}: \alpha, \beta\in A\}$.
\end{proof}

\noindent
{\em Proof of Theorem $1.2$.}
We consider the family of all
maps
$\{h_{\alpha}\colon Y_{\alpha}\to\uin^{\tau}\}_{
\alpha\in\Lambda}$ such that each $Y_{\alpha}$ is a closed subset
of $\uin^{\tau}$  with $Y_{\alpha}\in\mathcal{K}$.  Let $Y$ be the
disjoint sum of all  $Y_{\alpha}$ and the map $h\colon Y\to
\uin^{\tau}$ coincides with  $h_{\alpha}$ on every $Y_{\alpha}$. We extend 
$h$ to a map
$\overline{h}\colon\beta Y\to\uin^{\tau}$. Since
$\beta Y\in\mathcal K$, by Proposition 2.4,
there exists a compactum $X$ of weight $\leq\tau$ and maps 
$p\colon\beta Y\to X$ and $f\colon X\to\uin^{\tau}$ such that 
$X\in\mathcal K$ and $f\circ p=\overline{h}$.

Let us show that $f$ is $\mathcal K$-invertible. Take a space $Z\in\mathcal K$ 
and a map
$g\colon Z\to\uin^{\tau}$. Considering  $\beta
Z$  and the  extension $\overline{g}\colon\beta Z\to\uin^{\tau}$
of $g$, we can assume that  $Z$ is compact.  We also can
assume that the weight of $Z$ is $\leq\tau$ ( otherwise we apply
again Proposition 2.4 to find a compact
space $T\in\mathcal K$ of weight $\leq\tau$ and maps $g_1\colon Z\to T$ and
$g_2\colon T\to\uin^{\tau}$ with $g_2\circ
g_1=g$, and  then consider the space $T$ and the map $g_2$
instead, respectively, of $Z$ and $g$).  Therefore, without loss of 
generality, we can assume that  $Z$ is a closed subset of
$\uin^{\tau}$. According to the definition of $Y$ and the map $h$,
there is an index $\alpha\in\Lambda$ such that $Z=Y_{\alpha}$ and
$g=h_{\alpha}$.  The restriction $p|Z\colon Z\to X$ is a
lifting of $g$, i.e. $f\circ (p|Z)=g$.


\section{$C^{\ast}$-algebras with the approximate $n$-th root property }

In this Section we investigate the behavior of the classes $\mathcal{AP}(n)$, $\mathcal{AP}_1(n)$ and $\mathcal{HP}(n)$
with respect to direct systems and then use the result to prove the existence 
of universal elements in the classes $\displaystyle\mathcal{AP}(n)_s$, 
$\mathcal{AP}_1(n)_s$ and $\mathcal{HP}(n)_s$.

When we refer to a unital $C^{\ast}$-subalgebra of a
unital $C^{\ast}$-algebra we always assume that the
inclusion is a unital $\ast$-homomorphism.
The product in the category of (unital) $C^{\ast}$-algebras, i.e. the 
$\ell^{\infty}$-direct sum, is denoted by $\prod\{ X_{t} \colon t \in T \}$.
For a given set $Y$ and a cardinal number $\tau$, the symbol
$\exp_{\tau}Y$ denotes the partially ordered (by inclusion) set of
all subsets of $Y$ of cardinality not exceeding $\tau$.

Recall that a direct system
${\mathcal S} = \{ X_{\alpha}, i_{\alpha}^{\beta}, A\}$ of
unital $C^{\ast}$-algebras consists of a partially ordered directed
indexing set $A$,
unital $C^{\ast}$-algebras $X_{\alpha}$, $\alpha \in A$, and
unital $\ast$-homomorphisms
$i_{\alpha}^{\beta} \colon X_{\alpha} \to X_{\beta}$,
defined for each
pair of indexes $\alpha ,\beta \in A$ with $\alpha \leq \beta$,
and satisfying
the condition $i_{\alpha}^{\gamma} =
i_{\beta}^{\gamma}\circ i_{\alpha}^{\beta}$ for
each triple of indexes
$\alpha ,\beta ,\gamma \in A$ with $\alpha \leq \beta \leq \gamma$.
The (inductive) limit of the above direct
system is a unital $C^{\ast}$-algebra which is denoted by
$\varinjlim{\mathcal S}$. For each $\alpha \in A$ there
exists a unital $\ast$-homomorphism
$i_{\alpha} \colon X_{\alpha} \to \varinjlim{\mathcal S}$
which will be called
the $\alpha$-th limit homomorphism of $\mathcal S$.

If $A^{\prime}$ is a
directed subset of the indexing set $A$, then the subsystem
$\{ X_{\alpha}, i_{\alpha}^{\beta}, A^{\prime}\}$ of
${\mathcal S}$ is denoted ${\mathcal S}|A^{\prime}$.

Let $\tau \geq \omega$
be a cardinal number. A direct system
${\mathcal S} = \{ X_{\alpha}, i_{\alpha}^{\beta}, A\}$ of
unital $C^{\ast}$-algebras $X_{\alpha}$
and unital $\ast$-ho\-mo\-morp\-hisms
$i_{\alpha}^{\beta} \colon X_{\alpha} \to X_{\beta}$ is called a {\em direct
$C_{\tau}^{\ast}$-system} \cite{ch2} if the
following conditions are satisfied:
\begin{itemize}
\item[(a)]
$A$ is a $\tau$-complete set, that is, for each chain $C$ of elements of the 
directed set $A$ with 
${\mid}C{\mid} \leq \tau$, there exists an element
$\sup C$ in $A$. See \cite{ch1} for details.
\item[(b)]
The density $d(X_{\alpha})$ of $X_{\alpha}$ is at most $\tau$, for each 
$\alpha \in A$.
\item[(c)]
The $\alpha$-th limit homomorphism
$i_{\alpha} \colon X_{\alpha} \to \varinjlim{\mathcal S}$ is
an injective $\ast$-ho\-mo\-mor\-phism for each $\alpha \in A$.
\item[(d)]
If $B = \{ \alpha_{t} \colon t \in T\} $ is a chain of
elements of $A$ with
$|T| \leq \tau$ and $\alpha = \sup B$, then the limit
homomorphism
$\varinjlim\{ i_{\alpha_{t}}^{\alpha} \colon t \in T\}
\colon \varinjlim\left({\mathcal S}|B\right)
\to X_{\alpha}$ is an isomorphism.
\end{itemize}

\begin{pro}[Proposition 3.2, \cite{ch2}]
Let $\tau$ be an infinite cardinal number. Every unital
$C^{\ast}$-algebra $X$
can be represented as the limit of a direct
$C_{\tau}^{\ast}$-system
${\mathcal S}_{X} = \{ X_{\alpha}, i_{\alpha}^{\beta},
A \}$ where the index set $A = \exp_{\tau}Y$ for some (any)
dense subset $Y$ of $X$ with $|Y| = d(X)$.
\end{pro}

\begin{lem}[Lemma 3.3, \cite{ch2}]
If ${\mathcal S}_{X} = \{ X_{\alpha}, i_{\alpha}^{\beta}, A\}$
is a direct $C_{\tau}^{\ast}$-system, then
\[ \varinjlim{\mathcal S}_{X} = \bigcup\{
i_{\alpha}(X_{\alpha}) \colon \alpha \in A\} .\]
\end{lem}

The next proposition is a non-commutative version of Corollary 2.3

\begin{pro} Let $\mathcal K$ be one of the classes
$\mathcal{AP}(n)$, $\mathcal{AP}_1(n)$ and $\mathcal{HP}(n)$.
If $X$ is the limit of a direct system  
${\mathcal S} = \{ X_{\alpha}, i_{\alpha}^{\beta}, A\}$ consisting of unital 
$C^{\ast}$-algebras and unital $\ast$-inclusions with
$X_{\alpha}\in\mathcal K$ for each $\alpha$, then $X\in\mathcal K$.
\end{pro}

\begin{proof}  We consider first the case $\mathcal K=\mathcal{AP}(n)$. Let 
$a\in X$ with $\|a\|\leq 1$ and $\varepsilon>0$.
Since $\bigcup\{X_{\alpha} \colon \alpha \in A\}$ is dense in $X$ (we identify 
each $i_{\alpha}(X_{\alpha})$ with $X_{\alpha}$), there exist $\alpha$ and
$y\in X_{\alpha}$ with $\displaystyle\| a-y\|<\frac{\varepsilon}{4}$. Then, 
$\displaystyle\|y\|<\|a\|+\frac{\varepsilon}{4}\leq 1+\frac{\varepsilon}{4}$, 
so $\displaystyle\|\frac{y}{1+\varepsilon/4}\|<1$. Since 
$\displaystyle X_\alpha\in\mathcal{AP}(n)$, there is $b\in X_\alpha$ with 
$\displaystyle\|\frac{y}{1+\varepsilon/4}-b^n\|<\frac{\varepsilon}{2}$ and 
$\|b\|\leq 1$. Then 
$\displaystyle\|a-b^n\|\leq\|a-\frac{y}{1+\varepsilon/4}\|+\|\frac{y}{1+\varepsilon/4}-b^n\|<\varepsilon$. 
Hence, $X\in\mathcal{AP}(n)$. The above arguments work also for the class 
$\mathcal{HP}(n)$ because of the fact that the set of invertible elements of 
a $C^{\ast}$-algebra is open. Indeed, for an invertible element $a$ of $X$, 
the above fact allows us to choose $y$ in the above argument as an invertible 
element of $X$. Consequently, 
$\displaystyle\frac{y}{1+\varepsilon/4}$ is invertible in $X_\alpha$ and, 
since $X_\alpha\in\mathcal{HP}(n)$, there is $b\in X_\alpha$ with the required 
properties. Because the limit of any direct system consisting of 
$C^{\ast}$-algebras with bounded rank $\leq 1$ has a bounded rank $\leq 1$ 
\cite[Proposition 4.1]{cv}, the above proof remains valid for the class 
$\mathcal{AP}_1(n)$.
\end{proof}

As in the commutative case (see Proposition 2.6), we can establish a decomposition theorem for the classes $\mathcal{AP}(n)$, $\mathcal{AP}_1(n)$ and $\mathcal{HP}(n)$.

\begin{pro} Let $\mathcal K$ be one of the classes $\mathcal{AP}(n)$, 
$\mathcal{AP}_1(n)$ and $\mathcal{HP}(n)$.
The following conditions
are equivalent for any unital $C^{\ast}$-algebra $X$:
\begin{enumerate}
\item[(1)]
$X \in\mathcal K$.
\item[(2)]
$X$ can be represented as the direct limit of a direct
$C_{\omega}^{\ast}$-system
$\{ X_{\alpha}, i_{\alpha}^{\beta}, A\}$ satisfying the following properties:
\begin{itemize}
\item[(a)]
The indexing set $A$ is cofinal and $\omega$-closed in the
$\omega$-complete set $\exp_{\omega}Y$
for some (any) dense subset $Y$ of $X$ such that $|Y| = d(X)$.
\item[(b)]
$X_{\alpha}$ is a (separable) $C^{\ast}$-subalgebra of $X$ with
$X_{\alpha}\in\mathcal K$, $\alpha \in A$.
\end{itemize}
\end{enumerate}
\end{pro}
\begin{proof}
A similar statement holds for the class of all $C^{\ast}$-algebras of bounded 
rank $\leq n$ (see \cite[Proposition 4.2]{cv}). So, it suffices to consider 
the classes  $\mathcal{AP}(n)$ and $\mathcal{HP}(n)$. We suppose 
$\mathcal K=\mathcal{AP}(n)$.
The implication $(2) \Longrightarrow (1)$ follows from Proposition 3.3.

In order to prove the implication $(1) \Longrightarrow (2)$
we first consider a direct $C_{\omega}^{\ast}$-system
${\mathcal S}_{X} = \{ X_{\alpha}, i_{\alpha}^{\beta}, \Lambda\}$
with the properties indicated in Proposition 3.1 . Each $X_{\alpha}$ is 
identified with $i_{\alpha}(X_{\alpha})$. We next introduce the following 
relation $L \subseteq A^{2}$:\\
$(\alpha ,\beta ) \in \Lambda^{2}$ if and only if 
$\alpha \leq \beta$ and for each $x\in X_{\alpha}$ with $\|x\|\leq 1$ and
each $\varepsilon > 0$
there exists $y\in X_{\beta}$ such that 
$\|y\|\leq 1$ and $\|x-y^n\|<\varepsilon$.

Let us show that $L$ satisfies the following conditions:\\
(i) for every $\alpha\in \Lambda$ there exists $\beta\in \Lambda$ with 
$(\alpha,\beta)\in L$:\\ 
(ii) If $(\alpha,\beta)\in L$ and $\beta\leq\gamma$, then 
$(\alpha,\gamma)\in L$ ;\\ 
(iii) if $\{\alpha_k\}$ is a 
chain in $\Lambda$ with each $(\alpha_k,\beta)\in L$, then 
$(\alpha,\beta)\in L$, where $\alpha=\sup\{\alpha_k\}$.\\

To verify (i), we take $\alpha\in \Lambda$ and a countable set 
$M\subset X_{\alpha}$ which is dense in the unit ball 
$B_\alpha=\{x\in X_\alpha:\|x\|\leq 1\}$.  Since $X \in \mathcal{AP}(n)$, 
for each $x\in M$ and each $r\in Q^+$, we may take (and fix) 
$y(x,r)\in X$ with $\|x-y(x,r)^n\|<r$ and $\|y(x,r)\|\leq 1$. 
By Lemma 3.2, every $y(x,r)$ belongs to some $X_{\alpha (x,r)}$. Since 
$\Lambda$ is $\omega$-complete, according to \cite[Corollary 1.1.28]{ch1}, 
there exists $\beta\in \Lambda$ such that 
$\beta\geq\alpha$ and $\beta\geq\alpha(x,r)$ for each $x\in M$ and $r\in Q^+$. 
Then, $X_{\beta}$ contains all $y(x,r)$ and $(\alpha,\beta)\in L$. 
Condition (ii) follows directly because $\beta\leq\gamma$ implies 
$X_\beta\subset X_\gamma$. Let us establish condition (iii). If $\alpha$ is 
the supremum of the countable chain $\{\alpha_k\}$, then $X_\alpha$ is the 
direct limit of the direct system generated by the $C^{\ast}$-subalgebras 
$X_{\alpha_k}$, $k=1,2,..$, and the corresponding inclusion homomorphisms. 
This fact and $(\alpha_k,\beta)\in L$ for all $k$ yield $(\alpha,\beta)\in L$.

Since $L$ satisfies the conditions (i)-(iii), we can apply 
\cite[Proposition 1.1.29]{ch1} to conclude that the set 
$A=\{\alpha\in \Lambda:(\alpha,\alpha)\in L\}$ is cofinal and $\omega$-closed 
in $\Lambda$. Note that $(\alpha,\alpha)\in L$ precisely when 
$X_\alpha\in\mathcal{AP}(n)$. Therefore, we obtain a direct 
$C_{\omega}^{\ast}$-system 
${\mathcal S}_{X}^{\prime} = \{ X_{\alpha}, i_{\alpha}^{\beta}, A\}$ 
consisting of $C^{\ast}$-subalgebras $X_{\alpha} \in \mathcal{AP}(n)$ of 
$X$. Clearly 
$\varinjlim{\mathcal S}_{X}^{\prime} = X$. This completes the proof for the 
class $\mathcal{AP}(n)$. The case $\mathcal K=\mathcal{AP}(n)$ is similar.
\end{proof}

\noindent
{\em Proof of Theorem $1.4$.}
Let
${\mathcal B} =
\{ f_{t} \colon C^{\ast}\left( {\mathbb F}_{\infty}\right)
\to X_{t} \colon t \in T\}$ denote the set of all unital
$\ast$-ho\-mo\-mor\-phisms on
$C^{\ast}\left( {\mathbb F}_{\infty}\right)$
such that $X_{t}\in\mathcal K$. We claim that the product
$\prod\{ X_{t} \colon t \in T\}$ belongs to $\mathcal K$. This is obviously 
true if $\mathcal K$ is either $\mathcal{AP}(n)$ or $\mathcal{HP}(n)$. Since 
the bounded rank of this product is $\leq 1$ provided each $X_t$ is of 
bounded rank $\leq 1$ \cite[Proposition 3.16]{cv}, the claim holds for the 
class  $\mathcal{AP}_1(n)$ as well.
The $\ast$-homomorphisms $f_{t}$, $t \in T$, define the
unital $\ast$-homomorphism
$f \colon C^{\ast}\left( {\mathbb F}_{\infty}\right)
\to \prod\{ X_{t} \colon t \in T\}$ such that 
$\pi_{t} \circ f = f_{t}$ for each $t \in T$, where
$\pi_{t} \colon \prod\{ X_{t} \colon t \in T\} \to X_{t}$
denotes the canonical projection $\ast$-homomorphism onto $X_t$.
By Proposition 3.4,
$\prod\{ X_{t} \colon t \in T\}$ can be
represented as the limit of the $C^{\ast}_{\omega}$-system
${\mathcal S} = \{ C_{\alpha}, i_{\alpha}^{\beta}, A\}$ such
that $C_{\alpha}$ is a separable unital $C^{\ast}$-algebra
with $C_{\alpha}\in\mathcal K$ for each $\alpha \in A$. Suppressing 
the injective unital
$\ast$-homomorphisms $i_{\alpha}^{\beta} \colon C_{\alpha} \to C_{\beta}$, 
we may assume, for notational simplicity, that $C_{\alpha}$'s are unital
$C^{\ast}$-subalgebras of $\prod\{ X_{t} \colon t \in T\}$.
Let $\{ a_{k} \colon k \in \omega\}$ be a countable dense subset of
$C^{\ast}\left( {\mathbb F}_{\infty}\right)$. By Lemma 3.2,
for each $k \in \omega$ there exists an index $\alpha_{k} \in A$ such that
$f(a_{k}) \in C_{\alpha_{k}}$. Since $A$ is $\omega$-complete, there exists 
an index
$\alpha_0 \in A$ such that $\alpha_0 \geq \alpha_{k}$ for each $k \in \omega$.
Then
$f(a_{k}) \in C_{\alpha_{k}} \subseteq C_{\alpha_{0}}$ for each $k \in \omega$.
This observation coupled with the continuity of $f$ guarantees that
$f\left( C^{\ast}\left( {\mathbb F}_{\infty}\right)\right) =
f\left(\operatorname{cl}\left \{ a_{k} \colon k \in \omega\right\}\right)
\subseteq \operatorname{cl}\left\{ f\left(\{ a_{k} \colon k \in \omega\}
\right)\right\} \subseteq \operatorname{cl}C_{\alpha_{0}} = C_{\alpha_{0}}$.

Let $Z_{\mathcal K}= C_{\alpha_{0}}$ and define the unital $\ast$-homomorphism
$p \colon C^{\ast}\left( {\mathbb F}_{\infty}\right) \to Z_{\mathcal K}$ as 
$f$, regarded as a homomorphism of
$C^{\ast}\left( {\mathbb F}_{\infty}\right)$ into $Z_{\mathcal K}$.
Note that $f = i \circ p$, where
$i \colon Z_{\mathcal K}= C_{\alpha_{0}}
\hookrightarrow \prod\{ X_{t} \colon t \in T\}$ stands
for the inclusion.

By construction, we see $Z_{\mathcal K}\in \mathcal K$. Let us show that
$\displaystyle p \colon C^{\ast}\left( {\mathbb F}_{\infty}\right) \to 
Z_{\mathcal K}$
is $\mathcal K$-invertible.
For a given unital $\ast$-homo\-mor\-phism
$g \colon C^{\ast}\left( {\mathbb F}_{\infty}\right) \to X$, where $X$ is a
separable unital $C^{\ast}$-algebra with $X\in\mathcal K$, we
need to establish the existence of a unital $\ast$-homomorphism
$h \colon Z_{\mathcal K}\to X$ such that $g = h\circ p$. Indeed, by definition
of the set ${\mathcal B}$, we conclude that 
$g = f_{t} \colon C^{\ast}\left( {\mathbb F}_{\infty}\right) \to X_{t}= X$ 
for some index $t \in T$.
Observe that $g = f_{t} = \pi_{t} \circ f = \pi_{t} \circ i \circ p$.
This allows us to define the required unital $\ast$-homo\-mor\-phism
$h \colon Z_{\mathcal K}\to X$ as the composition $h = \pi_{t} \circ i$. Hence,
$p$ is $\mathcal K$-invertible.

\section{Example}

In this section, we show that a construction due to B. Cole 
(cf.\cite[Chap.3, section 19]{stout}) and M. Karahanjan 
\cite[Thoerem 5]{kar} yields a square root closed compatum $X$ 
such that $\displaystyle\check H ^1 (X; \mathbb Z)$ is infinitely generated.
In the sequel, we shall omit the coefficient group $\mathbb Z$. 
We will need the following theorem which is a consequence of 
\cite[Theorem 3.2]{gt}. 

\begin{thm}
Let $f\colon X\to Y$ be an open surjective map between compacta. 
Then $f^{*} \colon \check H ^1 (Y)\to\check H ^1 (X)$ is a
monomorphism.
\end{thm}

Now we outline the construction due to B.~Cole. This is based on
the exposition in \cite[Chapter 3, \S 19, p.194-197]{stout}. Let
$X$ be a compactum and define

$$S_X = \{(x, (z_f)_{f\in C(X)}) \colon f(x) = z^2_f
\,\,\,\,\mbox{ for each } f\in
C(X)\}\subset X\times \mathbb{C} ^{C(X)}$$

\noindent Note that $S_X$ is a closed subset of
$X\times\prod\{f(X)| f\in C(X)\}$ and hence is a compactum. Also, it is easy 
to see that $S_X$ is a pull-back in the
following diagram:

$$
\begin{diagram}
\node{S_X}\arrow{e}\arrow{s}\node{\mathbb{C} ^{C(X)}}\arrow{s,r}{S}\\
\node{X}\arrow{e,t}{F}\node{\mathbb{C} ^{C(X)}}
\end{diagram}
$$

\noindent 
where $F\colon X \to \mathbb{C}$ is defined by $F(x)=
(f(x))_{f\in C(X)} (x \in X),$ and 
$S\colon \mathbb{C} ^{C(X)}\to \mathbb{C}^{C(X)}$ is defined by 
$S((z_f)_{f\in C(X)})= (z^2_f)_{f\in C(X)}$.

 Let $\pi\colon S_X\to X$ be the map defined by $\pi [(x,
(z_f)_{f\in C(X)})]=x$ for all $x\in X$. Then $\pi$ is an open map
with zero-dimensional fibers. The critical property of $S_X$ and
$\pi$ is the following:\\
$(*)$ for any $f\in C(X)$ there exists $g\in C(X)$ such that $f\circ\pi = g^2$.
\\
Indeed, define $g\colon S_X\to\mathbb{C}$ by 
$g [(x, (z_f)_{f\in C(X)})] =z_f$. 

Note that $(*)$ implies:\\
$(**)$ $\pi$ is invertible with respect to the class of square root closed 
compacta.

Starting with a compactum $X_0$, by transfinite
induction we define an inverse spectrum 
$\{ X_{\alpha}, \pi^{\beta}_{\alpha}\colon X_{\beta}\to X_{\alpha} :
\alpha \leq \beta <\omega _1\}$ as follows. If $\beta =\alpha +1$ then 
$X_{\beta} = S_{X_{\alpha}}$ and $\pi _{\alpha} =\pi\colon X_{\beta} =
S_{X_{\alpha}} \to  X_{\alpha}$ is the map defined above. If
$\beta$ is a limit ordinal, then $\displaystyle X_{\beta}
=\lim_{\longleftarrow}(X_{\alpha}, \pi^{\gamma}_{\alpha}\colon
X_{\gamma}\to X_{\alpha}: \alpha \leq \gamma < \beta)$ and, for 
$\alpha < \beta$, let 
$\displaystyle\pi ^{\beta}_{\alpha} =\lim_{\longleftarrow}
(\pi^{\gamma}_{\alpha}\colon X_{\gamma}\to X_{\alpha}: \gamma < \beta).$

We let $X_{\Omega} = \displaystyle\lim_{\longleftarrow}
X_{\alpha}$. The $\alpha$-th limit projection is denoted by 
$\pi _{\alpha}:X_{\Omega} \to X_{\alpha}$. 
As the length of the above spectrum is $\omega _1$, the spectrum is 
factorizing in the 
sense that
each $f \in C(X_{\Omega})$ is represented as 
$f = f_{\alpha} \circ \pi_{\alpha}$ for some $\alpha < \omega _{1}$ and 
$f_{\alpha} \in C(X_{\alpha})$.
since its length is $\omega _1$. This implies that $C(X_{\Omega})$
is square root closed due to the property $(*)$.

In what follows, the unit disk in the complex plane 
$\{ z\in\mathbb{C} : |z|\le 1\}$ is denoted by $\Delta$.
 
\begin{thm}\label{infgen}  $C(\Delta_{\Omega})$ is square-root
closed, ${\rm dim} \Delta_{\Omega}\le 2$, 
$\displaystyle\check H ^1(\Delta _{\Omega})$ is infinitely generated and 
$2$-divisible.
\end{thm}

\noindent 
Notice that for each square root closed compactum $X$, 
$\displaystyle\check H ^1(X)$ is 2-divisible.  Hence, in view of the 
discussion above,
we need only to show that 
$\displaystyle\check H _1(\Delta_{\Omega})$ is infinitely generated. 
To show this, we need the following.

\begin{thm}\label{auxiliary}
$\check H ^1(S_{\Delta})$ is infinitely generated.
\end{thm}

\noindent 
Note that Theorem 4.2 
immediately follows from
Theorems 4.1 and 
Theorem 4.3. 
The proof of Theorem 4.3 
is divided into two parts.

\bigskip

\noindent Step 1. 
If  $\check H ^1(S_\Delta)$ is finitely generated
then $\check H ^1(S_\Delta)=0$.

\bigskip

\noindent Step 2. $\check H ^1(S_\Delta)\ne 0.$

\bigskip

\noindent Now we shall accomplish Steps 1 and 2.

\begin{pro}\label{retraction} Let $Y$ be a closed subspace of a compactum
$X$ such that there exists a retraction $r\colon X\to Y$.
Let also $i\colon Y\hookrightarrow X$ be the inclusion. Then there
exist an embedding $\overline{i} \colon S_Y\hookrightarrow S_X$ and a 
retraction
$\overline{r}\colon S_X\to S_Y$  such that the following diagram
is commutative.

$$
\begin{diagram}
 \node{S_Y}\arrow{e,t}{\displaystyle\overline{i}}\arrow{s,l}{\ds\pi
_{Y}}\node{S_X}\arrow{
e,t}{\ds\overline{r}}\arrow{s,l}{\ds\pi _X}\node{S_Y}\arrow{s,l}{\ds\pi _Y}\\
\node{Y}\arrow{e,t}{\ds i}\node{X}\arrow{e,t}{\ds r}\node{Y}
\end{diagram}
$$
\end{pro}

\begin{proof} Define $\overline{i}$ by

$$\overline{i}[ (y,(\eta _g)_{g\in C(Y)})]= (y, (\xi _f)_{f\in
C(X)})$$

\noindent where $\xi _f =\eta _{f|Y}$ for all $f\in C(X)$. Define
$\overline{r}$ by

$$\overline{r}[ (x,(\xi_f)_{f\in C(X)})]= (r(x), (\eta _g)_{g\in
C(X)})$$

\noindent where $\eta _g =\xi _{g\circ r}$ for all $g\in C(Y)$.
\end{proof}

\noindent Now we are ready to accomplish Step 1. Let $\Delta _m =
\{ z\in \mathbb{C} \colon |z|\le \frac{1}{m}\}\subset\Delta$. Let
$r_n\colon\Delta _n \to\Delta _{n+1}$ be the radial retraction and
$i_n\colon\Delta_{n+1}\hookrightarrow\Delta _n$ be the inclusion.
Consider the following sequence of commutative diagrams.

$$
\dgARROWLENGTH=2em
\begin{diagram}
\node{S_{\Delta_1}}\arrow{s,l}{\ds\pi_1} \node{S_{\Delta
_2}}\arrow{w,t}{\ds\overline{i}_1} \arrow{s,l}{\ds\pi _2}
\node{\dots} \arrow{w,t}{\ds\overline{i}_2}\node{S_{\Delta
_n}}\arrow{w} \arrow{s,l}{\ds\pi _n} \node{S_{\Delta
_{n+1}}}\arrow{w,t}{\ds\overline{i}_n} \node{\dots}\arrow{w}
\node{\lim _{\longleftarrow} S_{\Delta
_n}}\arrow{w}\arrow{s,l}{\ds\lim _{\longleftarrow}\pi_n = \pi
_{\infty}}
\\
\node{\Delta _1}  \node{\Delta _2}\arrow{w,t}{\ds i_1}
\node{\dots}\arrow{w,t}{\ds i_2} \node{\Delta _n}\arrow{w}
\node{\Delta _{n+1}}\arrow{w,t}{\ds i_n}  \node{\dots}\arrow{w}
\node{\{0\}}\arrow{w}
\end{diagram}
$$

\noindent It follows easily form the commutativity of the diagram
that $\lim _{\longleftarrow}S_{\Delta _n}$ is homeomorphic to the
inverse limit of the sequence

$$
\dgARROWLENGTH=2em
\begin{diagram}
\node{\pi^{-1}_1 (0)} \node{\pi^{-1}_2
(0)}\arrow{w,t}{\ds\overline{i}_1|} \node{\dots}\arrow{w}
\node{\pi^{-1}_{n} (0)}\arrow{w} \node{\pi^{-1}_{n+1}
(0)}\arrow{w,t}{\ds\overline{i}_n|} \node{\dots}\arrow{w}
\end{diagram}
$$

\noindent Since each fiber $\pi^{-1}_n (0)$ is $0$-dimensional, we
have ${\rm dim} \lim _{\longleftarrow} S_{\Delta _n} =0$. This
implies that $\ds \check H ^1 (\lim _{\longleftarrow} S_{\Delta
_n}) = \lim _{\longrightarrow}\check H^1( S_{\Delta _n}) =0$,
which is equivalent to the following observation.

\begin{pro}\label{triviallimit} For each 
$\alpha\in \check H ^1 (S _{\Delta _1}) = \check H
^1 (S _{\Delta })$, there exists an $n$ such that $(\overline{i}
_1\circ\dots \circ\overline{i}_n)^* (\alpha )=0$.
\end{pro}

\noindent Let $A_n$ be the annulus defined by $A_n
=\{z\in\mathbb{C} | \frac{1}{m+1}\le |z|\le\frac{1}{m}\}$, so that
$\Delta _n = \{0\}\cup (\cup\{ A_j |j\ge n\})$. Let $h\colon
\Delta =\Delta _1\to\Delta _2$ be the homeomorphism which maps
$A_j$ to $A _{j+1}$ ($j\ge 1$) by ``radial homeomorphisms" and
such that $h(0) =0$. Then the following diagram is commutative

$$
\begin{diagram}
\node{\Delta _n}\arrow{e,t}{\ds h|} \node{\Delta _{n+1}}\\
\node{\Delta _{n+1}}\arrow{n,l}{\ds i_n}  \arrow{e,t}{\ds h|}
\node{\Delta _{n+1}}\arrow{n,r}{\ds i_{n+1}}
\end{diagram}
$$

\noindent 
Define $h_n\colon S_{\Delta _{n}}\to S_{\Delta _{n+1}}$
by $h_n [ (x, (u_f)_{f\in C(\Delta _n )})] = (h(x), (v_g)_{g\in
C(\Delta _{n+1} )})$, where $v_g = u_{g\circ h}$, 
$g\in C(\Delta _{n+1})$. Note that $h_n$ is a homeomorphism.

\begin{pro}\label{commutative}  The following diagram is commutative.
$$
\begin{diagram}
\node{S_{\Delta _{n+1}}}\arrow{e,t}{\ds\overline{i}
_n}\arrow{s,l}{\ds h_{n+1}} \node{S
_{\Delta _n}}\arrow{s,r}{\ds h_{n}}\\
\node{S_{\Delta _{n+2}}}\arrow{e,t}{\ds\overline{i} _{n+1}}
\node{S _{\Delta _{n+1}}}
\end{diagram}
$$
\end{pro}
\begin{proof}
For each $(x_{n+1}, (z_f)_{f\in C(\Delta _{n+1} )})\in S_{\Delta
_{n+1}}$ we have
$$\overline{i} _n [(x_{n+1}, (z_f)_{f\in C(\Delta
_{n+1} )})] = (x_{n+1}, (u_f)_{f\in C(\Delta _{n} )})$$ where $u_f
= z_{f|\Delta _n} = z _{f\circ i_n}$, $f\in C(\Delta _n)$, and
$$h_n [(x_{n+1}, (u_f)_{f\in C(\Delta _{n} )})] = (h(x_{n+1}),
(v_f)_{f\in C(\Delta _{n+1} )})$$ where $v_f = u_{f\circ h} = z
_{(f\circ h)\circ i_n} = z _{f\circ ( h \circ i_n)}$. On the other
hand,
$$h_{n+1} [(x_{n+1}, (z_f)_{f\in C(\Delta _{n+1} )})] =
(h(x_{n+1}), (u_g)_{g\in C(\Delta _{n+1} )})$$ where $u_g =
z_{g\circ h}$, $g\in C(\Delta _{n+2})$, and
$$\overline{i} _{n+1}
[(h(x_{n+1}), (u_g)_{g\in C(\Delta _{n+1} )})] = (h(x_{n+1}),
(v_f)_{f\in C(\Delta _{n+1} )})$$
where $v_f = u_{f\circ i_{n+1}}
= z _{(f\circ i _{n+1})\circ h} = z _{f\circ (i_{n+1} \circ h})$.
Since $h\circ i_n = i_{n+1} \circ h$, we conclude that the diagram
is commutative.
\end{proof}

\noindent The above lemma provides a commutative diagram in
cohomologies:

\begin{center}
$$
\begin{diagram}
\node{\check H ^1(S_{\Delta _{n+1}})} \node{\check H ^1 (S
_{\Delta _n})} \arrow{w,t}{\ds\overline{i} _n ^*}
\\
\node{\check H ^1 (S_{\Delta _{n+2}})} \arrow{n,l}{\ds h_{n+1} ^*}
\node{\check H ^1 (S _{\Delta _{n+1}})} \arrow{n,r}{\ds h_{n} ^*}
\arrow{w,t}{\ds\overline{i} _{n+1} ^*}
\end{diagram}\leqno{(\dag)}
$$
\end{center}

\noindent 
Let 
$\phi = h^* _1\circ i ^* _1\colon \check H^1(S_{\Delta })\to 
\check H ^1 (S_{\Delta })$. Since $\overline
{r} _1\circ\overline{i} _1 ={\rm id}_{S_{\Delta }}$ we have
$\overline{i} _1^{\, *} \circ\overline {r} _1^*={\rm id}_{\check
H^1(S_{\Delta })}$ and hence $\phi$ is an epimorphism. We use
diagram $(\dag)$ to obtain the following diagram, in which all
vertical arrows are isomorphisms.

$$
\begin{diagram}
\node{\check H^1 (S_{\Delta})} \arrow{se,b}{\ds\phi}
\arrow{e,t}{\ds\overline{i}_1 ^{\, *}} \node{\check H^1 (S_{\Delta
_{2}})} \arrow{s,r}{h_1^*} \arrow{e,t}{\ds\overline{i}_2 ^{\, *}}
\node{\check H^1 (S_{\Delta _{3}})} \arrow{s,r}{h_2^*}
\arrow{e,t}{\ds\overline{i}_3 ^{\, *}} \node{\check H^1 (S_{\Delta
_{4}})} \arrow{s,r}{h_3^*} \arrow{e} \node{\dots}\\
\node{} \node{\check H^1 (S_{\Delta})} \arrow{se,b}{\ds\phi}
\arrow{e,t}{\ds\overline{i}_1 ^{\, *}} \node{\check H^1 (S_{\Delta
_{2}})} \arrow{s,r}{h_1^*} \arrow{e,t}{\ds\overline{i}_2 ^{\, *}}
\node{\check H^1 (S_{\Delta _{3}})} \arrow{s,r}{h_2^*} \arrow{e}
\node{\dots}
\\
\node{}\node{} \node{\check H^1 (S_{\Delta})}
\arrow{se,b}{\ds\phi} \arrow{e,t}{\ds\overline{i}_1 ^{\, *}}
\node{\check H^1 (S_{\Delta _{2}})} \arrow{s,r}{h_1^*}
\arrow{e} \node{\dots}\\
\node{}\node{}\node{} \node{\check H^1 (S_{\Delta})}
\arrow{se,b}{\ds\phi} \arrow{e}
\node{\dots}\\
\node{}\node{}\node{} \node{} \node{\vdots}
\end{diagram}
$$

\noindent 
The above diagram together with Proposition 4.5
imply that, for each $\alpha\in\check H^1 (S_{\Delta})$, 
there exists $n$ such that $\phi ^n (\alpha) =0$. If
$\check H ^1 (S _{\Delta})$ were finitely generated, we then would
have $\check H^1 (S _{\Delta})= 0$ because of the following
observation.

\begin{pro}
Let $A$ be a finitely generated Abelian group. If there exists an
epimorphism $f\colon A\to A$ such that for any $a\in A$ there
exists $n$ with $f^n (a) =0$, then $A$ is trivial.
\end{pro}

\begin{proof}
Note that $f\otimes 1 _{\mathbb{Q}}\colon A\otimes \mathbb{Q} \to 
A\otimes \mathbb{Q}$ is
an epimorphism of a vector space $A\otimes \mathbb{Q}$, which is
finite-dimensional over $\mathbb{Q}$. Hence $f\otimes
1_{\mathbb{Q}}$ is an isomorphism with the property in the hypothesis. 
This implies ${\rm rank} A =0$ and
therefore $A$ is a finite Abelian group. Then $f$ is an
isomorphism and therefore $A =0$.
\end{proof}

\noindent 
Thus Step 1 is completed and we proceed to Step 2.

\begin{pro}\label{open}
For a continuous function $f \in C(X)$, let \\ 
$S_f= \{ (x,z) : f(x) = z^2\,\, \mbox{ for each } x\in X\} \subset X \times 
\mathbb C$. 
Let also $\pi _f\colon S_f\to X$ be the projection. Then the natural map
$p_f\colon S_X\to S_f$, $(x, (z_g)_{g\in C(X)})\mapsto (x,z_f)$ is
open. Thus we have the following diagram.

$$
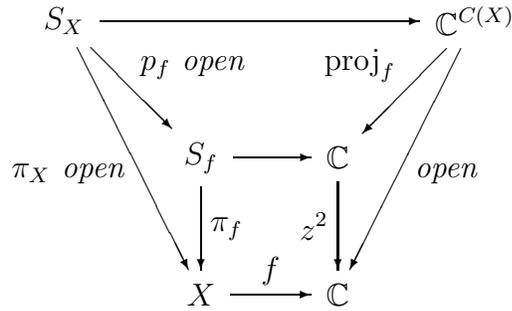
\begin{diagram}
\node{S_X}\arrow[3]{e}\arrow{se,t}{\ds p_f\mbox{
open}}\arrow{sse,b}{\ds \pi_X 
\mbox{ open}}\node[3]{\mathbb{C} ^{C(X)}}\arrow{sw,t}{\ds{\rm proj}_f}\arrow{ssw,b}{\mbox{open}}\\
\node{}\node{S_f}\arrow{e}\arrow{s,r}{\ds\pi _f}\node{\mathbb{C}}\arrow{s,l}{\ds z^2}\\
\node{}\node{X} \arrow{e,t}{\ds f}\node{\mathbb{C}}
\end{diagram}
$$
\end{pro}

\begin{proof}
Consider $g_1, g_2,\dots ,g_n\in C(X)$ and open subset $U_X\subset
X$,\\ $V_f, V_{g_1}, \dots ,V_{g_n}\subset\mathbb{C}$. It suffices
to show that

$$p_f [ (U_X\times V_f\times V_{g_1}\times \dots\times
V_{g_n}\times \prod _{g\ne g_1,\dots ,g_n,f} \mathbb{C})\cap
S_X]$$

\noindent 
is open in $S_f$. Take a point

$$(x, z_f, (z_{g_i})_{i=1}^n, (z_g)_{g\ne f,g_1,\dots ,g_n})\in
U_X\times V_f\times V_{g_1}\times \dots\times V_{g_n}\times \prod
_{g\ne g_1,\dots ,g_n,f} \mathbb{C}$$

\noindent 
and choose $\epsilon >0$ such that $B(z_f,\epsilon )
=\{w\in\mathbb{C}\colon |w-z_f|<\epsilon\}\subset V_f$ and
$B(z_{g_i},\epsilon)\subset V_{g_i}$ for all $i=1,2,\dots ,n$. Let
$a=f(x)$, $a_i = g_i (x)$, $i=1,2,\dots ,n$. There exists $\delta
>0$ such that if $|b-a|<\delta$ and $|b_i -a_i|<\delta$,
$i=1,\dots ,n$, then the equations 
\begin{eqnarray*}
z^2-b &=& 0 \\ 
z_i^2-b_i &=&0, i=1, \dots, n
\end{eqnarray*}
have solutions $z_b$ and $z_{b_i}$ respectively
such that $|z_b - z_f|<\epsilon$, $|z_{b_i} - z_{g_i}|<\epsilon$.
Choose a neighborhood $N$ of $x$ such that $|f(y) -f(x)|<\delta$
and $|g_i (y) - g_i (x)|<\delta$ for all $y\in N$ and $i=1,\dots,n$. We 
claim that

$$N\times B(z_f, \epsilon)\subset p_f [ (U_X\times V_f\times V_{g_1}\times 
\dots\times
V_{g_n}\times \prod _{g\ne g_1,\dots ,g_n,f} \mathbb{C})\cap
S_X]$$

\noindent Indeed, for each pint $(y,w)\in N\times B(z_f,
\epsilon)\subset N\times V_f$ we have $|g_i (y) - g_i
(x)|<\delta$, $i=1,2,\dots , n$ by the choice of $N$. Then we can
find $z_i\in B(z_{g_i},\epsilon )$ such that $z_i^2 = g_i (y)$.
Now for arbitrary choice of $z_g$, where $g\ne f,g_1,g_2,\dots
,g_n$ with $z^2_g = g(x)$, we have

$$(y, w, (z_i)_{i=1}^{n}, (z_g))\in U_X\times V_f\times V_{g_1}\times \dots\times
V_{g_n}\times \prod _{g\ne g_1,\dots ,g_n,f} \mathbb{C}$$

\noindent and $p_f [(y, w, (z_i)_{i=1}^{n}, (z_g))] = (y,w)$. This
proves the claim and hence completes the proof of the proposition.
\end{proof}

\noindent 
By Proposition 4.8 
and Theorem 4.1, 
the statement of the Step 2 follows from the next observation.

\begin{pro} There exists a mapping $f\colon \Delta \to \mathbb{C}$
such that $\check H^1 (S_f)\ne 0$.
\end{pro}

\begin{proof} Let $f(x,y) = (-2|x|+\sqrt{1-y^2}, y)$ for all
$(x,y)\in \Delta$. Then $S_f$ is homeomorphic to cylinder
$S^1\times I$.
\end{proof}

This completes the proof of Theorem 4.2.

\bigskip
The above construction is carried out word by word for disks of arbitrary 
dimensions.  In particular, applying the above to the one-dimensional disk 
$[-1,1]$, we have the following corollary which suggests that a topological 
characterization of general square root closed compacta could be rather 
different than the one for first-countable such compacta by \cite{hm} and 
\cite{mn}.

\begin{cor}
There exists an one-dimensional square root closed compactum $X$ with 
infinitely generated first {\v C}ech cohomology.
\end{cor}

For an infinite cardinal $\tau \geq \omega$, we consider 
$(\uin^{\tau})_{\Omega}$ 
and the limit projection 
$\pi _{\Omega}:(\uin^{\tau})_{\Omega} \to \uin^{\tau}$.  By the invertibility 
property (**) of $\pi:S_{X} \to X$ for arbitrary compactum $X$ and the 
standard spectral argument, it follows easily that 
$\pi _{\Omega}$ is also invertible with respect to the class of square root 
closed compacta. Hence we have 

\begin{pro}
The square root closed compactum $(\uin^{\tau})_{\Omega}$ contains every 
square root closed compactum of weight $\leq \tau$.
\end{pro}


\bigskip

\end{document}